\documentclass[11pt]{amsart}
\usepackage[pagebackref,colorlinks=true]{hyperref}
\usepackage{verbatim}
\usepackage{eucal,url,amssymb,stmaryrd,enumerate,amscd,}
\usepackage{amsfonts,amsmath,amsthm,amssymb,amscd,enumerate,eucal,url,stmaryrd}
\usepackage[margin=1in]{geometry}

\numberwithin{equation}{section}
\newtheorem{theorem}{Theorem}[section]

\newtheorem{proposition}[theorem]{Proposition}
\newtheorem{corollary}[theorem]{Corollary}
\newtheorem{definition}[theorem]{Definition}
\newtheorem{remark}[theorem]{Remark}
\newtheorem{example}[theorem]{Example}

\newtheorem{conv}[theorem]{Convention}

  \allowdisplaybreaks

\begin{document}

\begin{abstract}
Linear connections with torsion are important in the study of generalized Riemannian manifolds $(M,G=g+F)$, where the symmetric part $g$ of $G$ is a non-degenerate (0,2)-tensor and $F$ is the skew-symmetric part.
Some space-time models in theoretical physics are based on $(M,G=g+F)$, where $F$ is defined using an almost complex or almost contact metric structure.

In the paper, we first study more general models, where $F$ has constant rank and is based on weak metric structures (introduced by the first author and R.~Wolak), which generalize almost complex and almost contact metric structures.
We consider generalized metric connections (i.e., linear connections preserving $G$) with totally skew-symmetric torsion (0,3)-tensor.
For rank$(F)=\dim M$ and non-conformal tensor $A^2$, where $A$ is a skew-symmetric (1,1)-tensor adjoint to~$F$,
we apply weak almost Hermitian structures to fundamental results (by the second author and S. Ivanov) on generalized Riemannian manifolds and prove that the manifold is a weighted product of several nearly K\"ahler manifolds corresponding to eigen-distributions of $A^2$.
For rank$(F)<\dim M$ we apply weak $f$-structures and obtain splitting results for generalized Riemannian manifolds.


\end{abstract}

\title[Weak metric structures on generalized Riemannian manifolds]{Weak metric structures \\ on generalized Riemannian manifolds}
\date{\today}
\author{Vladimir Rovenski}
\address[Vladimir Rovenski]{Department of Mathematics, University of Haifa,
%Mount Carmel,
3498838 Haifa, Israel}
\email{vrovenski@univ.haifa.ac.il}
\author{Milan Zlatanovi\'c}
\address[Milan Zlatanovi\'c]{Department of Mathematics, Faculty of Science and Mathematics, University of Ni\v s, Vi\v segradska 33, 18000 Ni\v s, Serbia}
\email{zlatmilan@yahoo.com}

\date{\today}
\maketitle

\setcounter{tocdepth}{3}

%\tableofcontents

%%%\author{xxxx...}
%%%\address[xxxx...]{Univ...}
%%%\email{xxxx...}

\section{Introduction}

The {metric-affine~geometry} (founded by E.\,Cartan) generalizes Riemannian geometry: it uses a linear connection with torsion, $\nabla$, instead of the Levi-Civita connection of a Riemannian metric $g$, and appears in such context as homogeneous and almost Hermitian manifolds, Finsler and generalized geometries, e.g.~\cite{mikes}.
The~important distinguished case is {Riemann-Cartan manifolds}, where metric connections ($\nabla g =0$) are used. Riemann-Cartan spaces are central in gauge theory of gravity, where the torsion tensor is involved in the {Cartan spin-connection equation}.
Starting from 1950, A. Einstein used the real non-symmetric basic tensor $G$, sometimes called generalized Riemannian metric/manifold. In this theory the symmetric part $g$ of the basic tensor $G = g + F$ is related to gravitation, and the skew-symmetric one $F$ to electromagnetism.
Another concept of Generali\-zed Riemannian metrics/manifolds is part of Gene\-ralized Geometry introduced by N.~Hitchin and developed by
%his student
M.~Gualtieri in \cite{HG-2011}.
In our paper we do not consider this theory.

In this paper, we investigate linear connections on a generalized Riemannian manifold $(M,G)$ with a gene\-ral assumption, that the
symmetric part $g$ of $G$ is non-degenerate.
%We investigate connections on weak generalized %Riemannian manifolds, under the general assumption that %the symmetric part $g$ of $G$ is non-degenerate.
We extend examples of generalized Riemannian manifolds by weak almost (para) Hermitian manifolds,
weak (para) $f$-manifolds
and weak almost (para) contact manifolds,
%weak almost para-Hermitian, weak para-$f$-manifolds,
%and weak almost para-contact manifolds
defined by V. Rovenski and R. Wolak,~\cite{rov-108}.

{%\color{blue}
A \textit{weak almost Hermitian manifold} $(M, g, A, Q)$ is a Riemannian manifold $(M, g)$ of dimension $n\,(= 2m \ge 4)$ endowed with
non-singular endomorphisms: $A$ (skew-symmetric) and $Q$ (self-adjoint)
and the
%fundamental
2-form \(F\) such~that, see \cite{rov-137},
%the following conditions hold, see \cite{rov-137}:
\begin{equation}\label{WAH}
A^2 = -Q, \quad g(AX, AY) = g(QX, Y),
\quad F(X,Y)=g(AX,Y) .
\end{equation}
From $A^2=-Q$ we conclude that $A$ commutes with $Q$,
%\[
% $[A,Q]=0$;
%\]
thus
%\begin{align*}
 $F(X, QY) = F(QX, Y).$
%\end{align*}
 A \textit{weak almost contact metric manifold} $(M, A, Q, \eta,\xi, g)$ is a $(2m+1)$-dimensional Riemannian manifold equipped with a (1,1)-tensor $A$ of rank $2m$,
a non-singular (1,1)-tensor field $Q\ne{\rm Id}$,
a 1-form $\eta$ and a vector field $\xi$ dual to $\eta$ with respect to the metric $g$,\
$\eta(\xi)=1,\ \eta(X)=g(X,\xi)$,
satisfying
%the following compatibility conditions:
\begin{align}\label{acon}
A^2=-Q{+}\eta\otimes\xi, \
g(AX,AY)=g(QX,Y){-}\eta(X)\eta(Y), \
%F(X,Y)=g(AX,Y), \
A\xi=0, \
Q\xi=\xi.
\end{align}
In the case of \eqref{acon}, the skew-symmetric part $F$,
$F(X,Y)=g(AX,Y)$, of $G=g+F$ is degenerate,
$F(\xi,X)=0$, and has rank~$2m$.
%%%
Weak almost contact structure generalizes the classical almost contact structure, i.e., the~linear complex structure on the contact distribu\-tion is replaced by a nonsingular skew-symmetric (1,1)-tensor, hence the identity operator in the structural equations is replaced by a nonsingular self-adjoint (1,1)-tensor $Q$.
%allowed us to take a new look at the theory of contact manifolds.
Similarly we define a weak $f$-structure, see \eqref{Eq-f}, and the counterparts of the above:
weak para-Hermitian structure and weak para-$f$-structure, see Section~\ref{sec:3.3}.
}

The theory of weak metric structures, see \cite{rst-139}, is fully consistent with the skew-symmetric part $F$ of $G$; therefore it provides new tools for studying genera\-lized Riemannian manifolds.
We consider linear connections with torsion preserving the non-symmetric tensor $G$, so it means that it preserves its symmetric and skew-symmetric parts, $\nabla g=\nabla F=0$,
%{%\color{blue}
and consequently it preserves $Q$, $\nabla Q=0$,
where $Q$ is a self-adjoint endomorphism with respect to $g$ and $Q$ takes different values depending on the various types of manifolds mentioned~above.

In many geometric and physical contexts, such as string theory, gravity, and complex or almost complex geometry, skew-symmetric tensor $F$ plays a fundamental role  \cite{BSethi, FI, IZ1, IZ2, Iv, Str}.
Genera\-lized Riemanian manifolds studied here allow the natural incorporation of metric connections with torsion, particularly those preserving both $g$ and $F$, and provide a unified framework for exploring richer geometric structures beyond the scope of (pseudo-) Riemannian~geometry.

We pay special attention to the case when the torsion tensor of a linear
connection preserving the genera\-lized Riemannian metric $G$ is totally skew-symmetric w.r.t. the symmetric part $g$ of $G$ and satisfies the $A$-torsion or $Q$-torsion condition, depending on the structure.
The corresponding new Nijenhuis tensors na\-tu\-rally appear in weak almost (para) Hermitian, and weak almost (para) contact manifolds, whose total skew-symmetry depends on the $A$-torsion condition or the $Q$-torsion~condition.
Codazzi connection is defined on $(M, G=g+F)$,
%ge\-ne\-ralized Riemannian manifolds,
and an example of such a connection is given. Chern and Bismut connections are also studied in the context of weak Hermitian structures.

The paper is organized as follows.
Section~\ref{geomod} considers the geometric model. We introduce Codazzi connections and metric connections adapted to genera\-lized Riemannian manifolds, and discuss the torsion and Nijenhuis tensors, including their properties and significance.
Section~\ref{sec:03} presents applications of weak metric structures: weak almost Hermitian structures, weak $f$-structures,
%satisfying $f^3 + fQ = 0$,
 weak almost contact metric structures (Sections 3.1 and 3.2), and weak para-Hermitian and weak para-$f$-structures (Section 3.3)
% with $f^3 - fQ = 0$,
to generalized Riemannian manifolds.

\section{Linear Connections on Generalized Riemannian Manifolds}
\label{geomod}

\subsection{Geometrical Model}

The fundamental (0,2)-tensor $G$ in a non-symmetric (generalized) Riemannian manifold $(M,G)$ is in general non-symmetric. It is decomposed in two parts, the symmetric part $g$ (called (pseudo-) Riemannian metric) and the skew-symmetric part $F$ (called fundamental 2-form), $G=g+F$, where
\begin{equation*}%\label{metric}
g(X,Y)=\frac12\big[G(X,Y)+G(Y,X)\big], \qquad F(X,Y)=\frac12\big[G(X,Y)-G(Y,X)\big].
 \end{equation*}
We assume that the symmetric part, $g$, is non-degenerate of arbitrary signature,
and the skew-symmetric part, $F$, has a constant rank, e.g., is non-degenerate.
Therefore, we obtain a well-defined (1,1)-tensor $A$ of constant rank determined by the following condition:
\begin{equation}\label{m1}
 g(AX,Y) = F(X,Y)\quad \text{for all } X, Y \in \mathfrak{X}_M.
\end{equation}
By the above, since $F$ is skew-symmetric, then the tensor $A$ is also skew-symmetric:
%\[
 $g(AX, Y) = -g(X, AY)$
 %\quad \text{
 for all  $X, Y \in \mathfrak{X}_M$.
%\]
We seek natural linear connections $\nabla$ on $M$ preserving the generali\-zed Riemannian metric $G,\ \nabla G=0$, with a torsion (1,2)-tensor $T(X,Y)=\nabla_XY-\nabla_YX-[X,Y]$.

\begin{conv}%\hfill\break\vspace{-15pt} %
\label{Conv-1}\rm
In the paper we shall use the capital Latin letters $X,Y,\ldots$
to denote smooth vector fields on
%a smooth manifold
$M$, which commute, $[X,Y]=0$. Hence, $T(X,Y)=\nabla_XY-\nabla_YX$.
\end{conv}

Let $\nabla^g$ denote the Levi-Civita connection corresponding to the symmetric non-degenerate (0,2)-tensor $g$. The J.-L. Koszul's formula reads
\begin{equation*}%\label{lcg}
2g(\nabla^g_XY,Z)=
%\frac12\big[
Xg(Y,Z)+Yg(X,Z)-Zg(Y,X).
\end{equation*}
Denote the torsion (0,3)-tensor with respect to $g$ by the same letter,
\begin{equation*}%\label{t}
 T(X,Y,Z):=g(T(X,Y),Z).
\end{equation*}

%\subsection{Linear Connections on Generalized Riemannian Manifolds}

A
%linear
connection $\nabla$ with a torsion tensor $T$ on
%a generalized Riemannian manifold
$(M,G=g+F)$ is completely determined by the torsion tensor
and the covariant derivative $\nabla g$,
% of the symmetric part $g$ of $G$,
see \cite[Theorem~2.2]{IZ1}.
%More precisely, the following theorem is valid.
%\begin{theorem}[see \cite{IZ1}]\label{gmain}
%Let $(M,G=g+F)$ be a generalized Riemannian manifold
%equipped with a linear connection $\nabla$ with a %torsion tensor $T$.
%Then $\nabla$ is uniquely determined through $\nabla^g$,
%$\nabla g$ and $T$ by the following formula:
%\begin{align}\label{cgencon}
%\notag 2g(\nabla_XY,Z)&=2g(\nabla^g_XY,Z)
%+T(X,Y,Z)+T(Z,X,Y)-T(Y,Z,X)\\
%&-(\nabla_X\,g)(Y,Z)-(\nabla_Y\,g)(Z,X)+(\nabla_Z\,g)(Y,X).
%\end{align}
%The covariant derivative $\nabla F$ of the skew-symmetric part $F$ of $G$ is given by
%\begin{align}\label{cconn2}
%\notag
% 2(\nabla_XF)(Y,Z)&=2(\nabla^g_XF)(Y,Z) %+T(X,Y,AZ)+T(Z,X,AY)\\
%&\notag
%+T(AZ,X,Y)+T(AZ,Y,X)+T(X,AY,Z) +T(Z,AY,X)\\
%\notag
%&+(\nabla_X\,g)(AY,Z)-(\nabla_X\,g)(Y,AZ)-(\nabla_Y\,g)(AZ,X)\\
%&+(\nabla_Z\,g)(AY,X)+(\nabla_{AZ}\,g)(Y,X)
%-(\nabla_{AY}\,g)(Z,X).
%\end{align}
%In particular, the exterior derivative $dF$ of $F$ satisfies
%\begin{align}\label{cconn1}
%\notag
%dF(X,Y,Z)=-T(X,Y,AZ)-T(Y,Z,AX)-T(Z,X,AY)\\
%&+(\nabla_XF)(Y,Z)+(\nabla_YF)(Z,X)+(\nabla_ZF)(X,Y).
%\end{align}
%Conversely, any three tensors $T,\,\nabla g$ and
%$\nabla F$ satisfying \eqref{cconn2} determine a unique %linear connection $\nabla$ via~\eqref{cgencon}.
%\end{theorem}
%
Using the vector fields defined in Convention~\ref{Conv-1}, the co-boundary formula for exterior derivative of a 2-form $F$ reduces to the following:
\begin{align*}%\label{E-3.3}
 \,d F(X,Y,Z) = X(F(Y,Z)) + Y(F(Z,X)) + Z(F(X,Y)).
% \\ &-F([X,Y],Z) - F([Z,X],Y) - F([Y,Z],X)
%\big].
\end{align*}
Here this formula is without the coefficient 3, unlike \cite{Blair}.

%{%\color{red}\% MZ We use the "coordinate system" as %defined in Convention 2.1, where any two vector fields %commute. Therefore, the last term will disappear. }

\begin{definition}\rm
Let $(M, G = g + F)$ be a generalized Riemannian manifold and $A \colon TM \to TM$ a skew-symmetric endomorphism of constant rank, see~\eqref{m1}.
A~linear connection $\nabla$ on $M$ is said to have the {\it A-torsion condition} if its torsion tensor $T$ satisfies the following:
\begin{equation}\label{Eq-A-T}
T(AX, Y) = T(X,AY)  \quad \text{for all } X, Y \in \mathfrak{X}_M .
\end{equation}
%%
%A linear connection \( \nabla \) is said to have {\it $Q$-torsion condition} if its torsion tensor \( T \) satisfies
%\begin{equation}\label{Eq-Q-T}
%T(QX, Y) = T(X, QY)  \quad \text{for all } X, Y \in %\mathfrak{X}_M,
%\end{equation}
%where \( Q: TM \to TM \) is an endomorphism that is self-%adjoint with respect to the symmetric part $g$ of metric %$G$, i.e. $g(QX, Y) = g(X, QY)$ for all  $X, Y \in %\mathfrak{X}_M$.
\end{definition}

\subsection{Codazzi connections}
 %on generalized Riemannian manifolds

A linear connection $\nabla$ on a generalized Riemannian manifold $(M,G=g+F)$ is called a \textit{Codazzi connection}
if the covariant derivative of the symmetric part of the metric tensor $G$ satisfies the following equation, e.g., \cite{Amari2016,PA-2020}:
\[
 (\nabla_Z\,g)(X,Y)=(\nabla_X\,g)(Z,Y) \quad \text{for all } X, Y, Z \in \mathfrak{X}_M .
\]
%for all vector fields \(X, Y, Z\).
Such a pair $(G,\nabla)$ is called a \textit{Codazzi structure} on $(M, G)$.
As a consequence, the tensor $\nabla g$ is completely symmetric in all three variables,  that is,
\begin{equation}\label{ggstat}
     (\nabla_Z\,g)(X,Y)= (\nabla_X\,g)(Z,Y) =(\nabla_Y\,g)(X,Z).
\end{equation}

This structure plays a fundamental role in information geometry and the theory of statistical connections, enabling a deep relationship between the geometric properties of the metric and the statistical properties of models.

\begin{theorem}%\label{T-gstat}
Let $(M,G=g+F)$ be a generalized Riemannian manifold equipped with
%$\nabla^g$ be the Levi-Civita connection of $g$.
%Let
a Codazzi connection $\nabla$ having a torsion tensor $T$.
% and denote the covariant derivative of the symmetric part $g$ of $G$ by $\nabla g$.
Then $\nabla$ is uniquely determined by the following formula:
\begin{align*}%\label{cstatcon}
 2g(\nabla_XY,Z)=2g(\nabla^g_XY,Z) +
 %\frac 12\,
 T(X,Y,Z)
 -
 %\frac12\,
 (\nabla_X\,g)(Y,Z).
\end{align*}
The covariant derivative $\nabla F$ (of the skew-symmetric part $F$) of $G$ is given by
\begin{align*}%\label{stat}
%\notag
2(\nabla_XF)(Y,Z)&=2(\nabla^g_XF)(Y,Z) +
%\frac12\big[
T(X,Y,AZ)+T(Z,X,AY)\\
&+
%\frac12\big[
T(AZ,X,Y)+T(AZ,Y,X)+T(X,AY,Z)+T(Z,AY,X)\\
& +
%\frac12\big[
(\nabla_Z\,g)(AY,X) -(\nabla_Y\,g)(AZ,X).
\end{align*}
\end{theorem}

\begin{proof}
Using \eqref{ggstat}, it holds
$$(\nabla_X\,g)(AY,Z)=(\nabla_{AY}\,g)(Z,X), \quad (\nabla_X\,g)(Y,AZ)=(\nabla_{AZ}\,g)(Y,X).$$
Then, by applying
%Theorem~\ref{gmain},
\cite[Theorem~2.2]{IZ1}
the result follows directly.
\end{proof}

%{%\color{blue}
\begin{example}\rm
In the sense of Eisenhart's definition, see~\cite{Eisen}, a generalized Riemannian manifold is a differentiable manifold equipped with a non-symmetric basic tensor
\( G(X,Y) = g(X,Y) + F(X,Y) \),
with a linear connection
explicitly defined by the equality
\begin{align}\label{EisCon}
\notag
 g(\nabla_X Y, Z) &= \frac{1}{2} \big[ X G(Y, Z) + Y G(Z, X) - Z G(Y, X) \big] \\
 & =g(\nabla^g_XY, Z) + \dfrac 12dF(X, Y, Z).
\end{align}
Using the definition of the connection~\eqref{EisCon},
we have the following:
\begin{align}\label{EinProp}
 g(\nabla_X Y, Z) + g(Y, \nabla_X Z)
 %&= \frac{1}{2} \big[ X G(Y,Z) + Y G(Z,X) - Z G(Y,X) %\nonumber \\
 %&\quad + X G(Z,Y) + Z G(Y,X) - Y G(Z,X) \big] \\
 %&= \frac{1}{2} \big[ X G(Y,Z) + X G(Z,Y) \big] %\nonumber \\
 &=Xg(Y,Z).
 %\nonumber
\end{align}
%Therefore,
%\[
%g(\nabla_X Y, Z) + g(Y, \nabla_X Z) = X g(Y,Z).
%\]
By using \eqref{EinProp}, the covariant derivative of the metric $g$ is the following:
\begin{align*}
%\notag
(\nabla_X g)(Y,Z) &=X(g(Y,Z))-g(\nabla_XY,Z) -g(Y,\nabla_XZ)\\
& =X(g(Y,Z))-X(g(Y,Z))=0.
\end{align*}
Thus, the metric $g$ is covariantly constant with respect to $\nabla$.
%, i.e.,  $\nabla g = 0$.
%%
The torsion (0,3)-tensor $T$ is totally skew-symmetric and is given by the equation $T(X, Y, Z) = dF(X, Y, Z)$.
Therefore, the connection on a generalized Riemannian manifold in the sense of Eisenhart's definition is a Codazzi connection.
\end{example}
%}

%\subsection{Skew-symmetric Codazzi-type structure.}

%A linear connection $\nabla$ on a generalized Riemannian manifold %$(M,G=g+F)$ is called \textit{skew-symmetric Codazzi-type}
%if the covariant derivative of the skew-symmetric part with respect to %Levi-Civita connection of $g$ satisfies the following equation
%\[
%(\nabla^g_Z F)(X,Y) = (\nabla^g_X F)(Y,Z)
%\]
%for all vector fields \(X, Y, Z\). Such a pair $(G,\nabla)$ is called a %\textit{skew-symmetric Codazzi-type structure} on $M$.

%We obtain that the covariant derivative with respect to Levi-Civita connection of the 2-form $F$, which is antisymmetric in its first two arguments, satisfies the cyclic identity, i.e.
%\begin{equation}\label{ssCodazi}(\nabla^g_Z F)(X,Y) = (\nabla^g_X F)(Y,Z)=%(\nabla^g_Y F)(Z,X).
%\end{equation}

\subsection{Generalized metric connections}
% on generalized Riemannian manifolds

Here, we consider
%investigate the existence of
generalized metric connections, i.e., linear connections $\nabla$ preserving the generalized Riemannian metric, $G$, on a generalized Riemannian manifold $(M,G=g+F)$. First, we recall the following result in \cite{IZ1}.

\begin{theorem}
%see \cite{IZ1}]
\label{main}
Let $(M,G=g+F)$ be a generalized Riemannian manifold.
%and $\nabla^g$ be the Levi-Civita connection of $g$.
\begin{enumerate}
\item[i)] A linear connection $\nabla$ on $M$ preserves the generalized Riemannian metric $G$ if and only if it  preserves its symmetric part $g$ and its skew-symmetric part $F$:
$%\begin{equation}\label{conn}
\nabla G=0 \Leftrightarrow \nabla g=\nabla F=0 \Leftrightarrow
\nabla g=\nabla A=0.
$ %\end{equation}
\item[ii)] If
%there exists
$\nabla$ is a generalized metric connection
%linear connection on $M$ preserving the generalized %Riemannian metric $G,\ \nabla G=0$,
with a torsion tensor $T$, then
%the following
%condition
%holds:
\begin{align}\label{conn2}
\notag
 2(\nabla^g_XF)(Y,Z)= & -T(X,Y,AZ)-T(Z,X,AY) -T(AZ,X,Y)\\
 &-T(AZ,Y,X)-T(X,AY,Z)-T(Z,AY,X);
\end{align}
in particular, the exterior derivative of $F$ satisfies the following equality:
\begin{align}\label{conn1}
 dF(X,Y,Z)=-T(X,Y,AZ)-T(Y,Z,AX)-T(Z,X,AY);
\end{align}
\indent
Conversely, if
%the condition
\eqref{conn2} is valid, then there exists a unique generalized metric connection $\nabla$ on $(M,G=g+F)$ determined by the torsion $T$ by
%the following formula:
\begin{equation}\label{gencon} 2g(\nabla_XY,Z)=2g(\nabla^g_XY,Z) +T(X,Y,Z)+T(Z,X,Y)
-T(Y,Z,X).
\end{equation}
\end{enumerate}
\end{theorem}

\begin{remark}\rm
According to i) of Theorem~\ref{main}, a generalized Riemannian manifold $(M,G=g+F)$ equipped with a generalized metric connection $\nabla$ can be seen as a Riemann-Cartan manifold $(M,g,\nabla)$ equipped with the fundamental 2-form~$F$.
The~formulas in Theorem~\ref{main} are especially meaningful under the assumption that the torsion tensor \(T\) satisfies the $A$-torsion condition (\ref{Eq-A-T}) with totally skew-symmetric property.
\end{remark}

\begin{proposition}
Let \((M, G = g + F)\) be a generalized Riemannian manifold with non-degenerate fundamental 2-form $F$.
Then a generalized metric connection $\nabla$ with a totally skew-symmetric torsion $(0,3)$-tensor $T$
%the linear connection \(\nabla\)
satisfies the $A$-torsion condition \eqref{Eq-A-T} if and only if $(\nabla^g_{AX} A)X=0$ for all $X\in\mathfrak{X}_M$.
%is valid.
\end{proposition}

%\begin{proof}
%Using~\eqref{conn2} and totally skew-symmetric torsion (0,3)-tensor, by replacing \(X\) with \(AX\) in the first term and \(Y\) with \(AY\)
%in the second term, we obtain
%\begin{align*}
%& (\nabla^g_{AX} F)(Y,Z) + (\nabla^g_{AY} F)(X,Z)= -\dfrac{1}{2} \big[ T(AX,Y,AZ)-T(X,AY,AZ) \big].
%\end{align*}
%Since \(A\) is a non-singular endomorphism, the above implies the assertion.
%\end{proof}

\begin{proof}
Using~\eqref{conn2} and totally skew-symmetric torsion (0,3)-tensor, by replacing \(X\) with \(AX\) in the first term and \(Y\) with \(AY\) in the second term, we obtain
\begin{align*}
& (\nabla^g_{AX} F)(Y,Z) + (\nabla^g_{AY} F)(X,Z)
%\\
%& = -\dfrac{1}{2} \big[ T(AX,Y,AZ) + T(Z,AX,AY) \big]
% - \dfrac{1}{2} \big[ T(AY,X,AZ) + T(Z,AY,AX) \big] \\
%&
= -\dfrac{1}{2} \big[ T(AX,Y,AZ)-T(X,AY,AZ) \big].
\end{align*}
{%\color{blue}
If $(\nabla^g_{AX} A)X=0$, then the left-hand side vanishes, hence
\[
T(AX,Y,AZ) = T(X,AY,AZ)\quad\text{i.e.} \quad AT(AX,Y)=AT(X,AY).
\]
Since $A$ is a non-singular endomorphism, this implies
%\[
$T(AX,Y) = T(X,AY)$.
%\]
Moreover, using the skew-symmetry again, we compute
\[
T(AX,Y,Z) = -T(AX,Z,Y) = T(X,AY,Z) = -T(X,Z,AY).
\]
By interchanging $Y$ and $Z$ in the second and fourth terms, it follows that
%\[
$T(AX,Y,Z) = T(X,Y,AZ)$,
%\]
and finally
\[
T(AX,Y,Z) =T(X,AY,Z)= T(X,Y,AZ),
\]
the above completes the proof.
}
\end{proof}

\begin{definition}\rm
Let $\nabla$ be a linear metric
%(non-symmetric)
connection on a generalized Riemannian mani\-fold and $\nabla^g$ the Levi-Civita connection of $g$.
The \emph{contorsion (or, difference)} (1,2)-tensor \( K \) is defined~by
%\begin{equation*}%\label{eq:contorsion-def}
 $K(X, Y) = \nabla_X Y - \nabla^g_X Y$.
%\end{equation*}
The \emph{contorsion} (0,3)-tensor is defined by
\begin{equation*}%\label{kt}
 K(X,Y,Z):=g(K(X,Y),Z).
\end{equation*}
\end{definition}

%The Levi-Civita connection is a distinguished point in affine space of all linear connections $\nabla$ on $M$ that is in one-to-one correspondence with the %linear space of all contorsion tensors $K$ on $(M,g)$.

By Theorem~\ref{main}, any generalized metric connection on a generalized Riemannian manifold is completely
determined by the torsion tensor~$T$.
If $\nabla$ is a generalized metric connection on
%a generalized Riemannian manifold,
$(M,G=g+F)$,
%preserves the generalized Riemannian metric $G$,
then from (\ref{gencon})
%, see also \eqref{cgencon} with $\nabla g=0$,
we~have
\begin{equation}\label{kktorz}
 2K(X,Y,Z)=
 %\frac12\big[
 T(X,Y,Z)+T(Z,X,Y)-T(Y,Z,X).
\end{equation}
In this case, the contorsion (0,3)-tensor $K(X,Y,Z)$ is totally skew-symmetric if and only if the torsion (0,3)-tensor $T(X,Y,Z)$ is totally skew-symmetric.

\begin{corollary}
 Let $\nabla$ be a linear connection with a totally skew-symmetric torsion $(0,3)$-tensor.
 Then the $A$-torsion condition \eqref{Eq-A-T} holds if and only if the following $K$-torsion condition holds:
\begin{equation*}%\label{Eq-K-T}
 K(AX, Y) = K(X,AY)=AK(X,Y)  \quad \text{for all } X, Y \in \mathfrak{X}_M .
\end{equation*}
\end{corollary}

%\begin{cor}
%Let $\nabla$ be a linear connection with a totally skew-symmetric torsion $(0,3)$-tensor. Then the $A$-torsion condition \eqref{Eq-A-T} holds if and only if %the following $K$-torsion condition holds:
%\begin{equation*}%\label{Eq-K-T}
% K(AX, Y) = K(X,AY)=AK(X,Y)  \quad \text{for all } X, Y \in \mathfrak{X}_M .
%\end{equation*}
%\end{cor}

{%\color{blue}
\begin{proof}
Since the torsion tensor is totally skew-symmetric, equation \eqref{kktorz} becomes
\begin{align}\label{ktt}
 K(X,Y,Z)
  &= \dfrac{1}{2} \big[ T(X,Y,Z) + T(Z,X,Y) - T(Y,Z,X) \big] \notag \\
  &= \dfrac{1}{2} \big[ T(X,Y,Z) + T(X,Y,Z) - T(X,Y,Z) \big] \notag \\
  &= \dfrac{1}{2} \, T(X,Y,Z).
\end{align}
If the $A$-torsion condition \eqref{Eq-A-T} for totally skew-symmetric torsion is true, then
\begin{equation}\label{Eq-T1}
T(AX, Y, Z) = -T(AX,Z,Y) = T(X,AY,Z) = -T(X,Z,AY).
\end{equation}
By interchanging $Y$ and $Z$ in the second and fourth terms of \eqref{Eq-T1}, we obtain
%\begin{equation*}%\label{Eq-T2}
$T(AX,Y,Z) = T(X,Y,AZ)$.
%\end{equation*}
Therefore,
\begin{equation}\label{Eq-T3}
T(AX,Y,Z) = T(X,AY,Z) = T(X,Y,AZ).
\end{equation}
The identity \eqref{Eq-T3} together with (\ref{ktt}) gives
\begin{equation*}%\label{Eq-K}
K(AX,Y,Z) = K(X,AY,Z) = K(X,Y,AZ),
\end{equation*}
which completes the proof.
\end{proof}
}

\subsection{The Nijenhuis tensor}

The Nijenhuis tensor $N_P$ of a (1,1)-tensor~$P$ on a smooth manifold $M$ is defined as, e.g.\cite{KN},
\begin{equation}\label{nuj}
 N_P(X,Y)=[PX,PY]+P^2[X,Y]-P[PX,Y]-P[X,PY].
\end{equation}
%see e.g.\cite{KN}.
The Nijenhuis tensor is skew-symmetric by definition and it plays a fundamental role in almost complex (resp. almost para-complex) geometry.
We denote the Nijenhuis (0,3)-tensor with respect to $g$ with the same letter,
\[
 N_P(X,Y,Z):=g(N_P(X,Y),Z).
\]
If $A^2=-{\rm Id}$ (resp. $A^2={\rm Id}$) then the celebrated
{%\color{blue}
Newlander-Nirenberg
}
theorem (see, e.g. \cite{KN}) shows that an almost complex structure  is integrable if and only if the Nijenhuis tensor $N_A$ vanishes.

Let $\nabla$ be a linear connection on $M$ preserving the generalized Riemannian metric $G,\ \nabla G=0$.
Then $\nabla$ preserves $g,F$ and
$A$:\ $\nabla g=\nabla F=\nabla A=0$.
Using the definition of the torsion tensor $T$ and the covariant derivative $\nabla A$, we express the Nijenhuis tensor in terms of $T$
%the torsion tensor
and $\nabla A$ as~follows:
\begin{align}\label{nuj1}
\notag
N_A(X,Y)&=(\nabla_{AX}A)Y-(\nabla_{AY}A)X-A(\nabla_{X}A)Y+A(\nabla_{Y}A)X\\
&-T(AX,AY)-A^2T(X,Y)+AT(AX,Y)+AT(X,AY).
\end{align}
Using $\nabla A=0$ in \eqref{nuj1} and applying \eqref{m1}, we get the following:
\begin{align}\label{nuj2}
\notag
 N_A(X,Y,Z) & = -T(AX,AY,Z) -T(X,Y,A^2Z) \\
 & -T(AX,Y,AZ) -T(X,AY,AZ).
\end{align}

%{%\color{blue}
%\begin{definition}\rm
%Let  \( Q: TM \to TM \) be a self-adjoint endomorphism with respect to the symmetric part $g$ of metric $G$, i.e., $g(QX, Y) = g(X, QY)$ for all $X, Y \in \mathfrak{X}_M$. The \emph{Nijenhuis tensor} of \( Q \) is the (1,2)-tensor defined by
%\begin{equation}\label{NujQ}
% N_Q(X,Y) = [QX, QY] + Q^2[X,Y] - Q[QX,Y] - Q[X,QY].
%\end{equation}
%\end{definition}

\subsection{The totally skew-symmetric torsion (0,3)-tensor}

%In this section we investigate the torsion tensor.
Pseudo-Riemannian manifolds equip\-ped with metric connections
%Linear connections preserving a pseudo-Riemannian
%metric and
having a totally skew-symmetric torsion (0,3)-tensor $T$, i.e. the following condition holds:
\begin{equation}\label{st}
 T(X,Y,Z)=-T(X,Z,Y),
\end{equation}
became very attractive in the last twenty years mainly due to the relations with supersymmetric string theories (see \cite{BSethi,GPap,Str} and references therein, for a mathematical treatment consult the nice overview \cite{Agr}). The main point is that the number of preserved supersymmetries is equal to the number of parallel spinors with respect to such a connection.

%In the following,
%We investigate when a generalized Riemannian manifold %admits a metric connection with totally skew-symmetric %torsion tensor.
%In the case, when a
For a generalized Riemannian manifold $(M,G=g+F)$ equipped with a generalized metric connection $\nabla$ that has a totally skew-symmetric torsion (0,3)-tensor, the following is true, see Theorem~2.5 in \cite{IZ1}:
%\begin{theorem}[see \cite{IZ1}]%\label{skew}
%Let $(M,G=g+F)$ be a generalized Riemannian manifold.
%If there exists a linear connection $\nabla$ preserving %the generalized Riemannian metric $G,\ \nabla G=0$,
%with totally skew-symmetric torsion tensor $T$,
%then the following condition holds:
\begin{equation}\label{ndf1}
 N_A(X,Y,AZ)+N_A(X,Z,AY)=dF(X,Y,A^2Z)+dF(X,Z,A^2Y).
\end{equation}
Moreover, the torsion (0,3)-tensor $T$ of $\nabla$ satisfies the following equalities:
\begin{equation}\label{tor1}
\begin{split}
T(AX,AY,Z)=-N_A(X,Y,Z)+dF(X,Y,AZ),\\
T(AX,Y,Z)=2(\nabla^g_XF)(Y,Z)-dF(X,Y,Z),
\end{split}
\end{equation}
and
%the torsion tensor of the connection $\nabla$
is determined by the following formula:
\begin{equation}\label{skct}
% 2g(\nabla_XY,Z)=2g(\nabla^g_XY,Z)+T(X,Y,Z).
 T(X,Y,Z)=2g(\nabla_XY,Z)-2g(\nabla^g_XY,Z).
\end{equation}
If, in addition, the fundamental 2-form
$F$ of the generalized Riemannian metric $G$ is closed, $dF=0$, then \eqref{ndf1} and \eqref{tor1} hold inserting $dF=0$.
%\end{theorem}

%\begin{lemma}
%Let $(M,G=g+F)$ be a generalized Riemannian manifold and
%$\nabla^g$ be the Levi-Civita connection of $g$. If there exists a
%linear connection $\nabla$ preserving the generalized Riemannian
%metric $G,\ \nabla G=0$ with totally skew-symmetric torsion $T$, then the %following two conditions are equivalent:
%\begin{enumerate}
%    \item $\nabla$ is skew-symmetric Codazzi-type.
%\item The  skew-symmetric torsion tensor $T$  satisfies
%    \begin{equation}\label{torsionA}
%       T(AX,Y,Z) = T(X,AY,Z) = T(X,Y,AZ).
%    \end{equation}
%\end{enumerate}
%\end{lemma}

%\noindent{\it Proof.} The equivalence follows by straightforward calculations from the second equation in in (\ref{tor1}).

%From (\ref{torsionA}) and (\ref{tor1}) immediately follows
%\begin{equation}\begin{split}
%T(AX,Y,Z)=-\dfrac13dF(X,Y,Z)\\
%N_A(X,Y,Z)=\dfrac 43 dF(X,Y,AZ)\\
%(\nabla^g_XF)(Y,Z)=\dfrac 13dF(X,Y,Z).\end{split}
%\end{equation}
%Consequently, the exterior derivative $dF$ of the 2-form $F$ satisfies %the following symmetry property:
%\begin{equation}\label{dfah}
%dF(X, Y, AZ) = dF(X, AY, Z) = dF(AX, Y, Z),
%\end{equation}

\section{New Applications of Weak Metric Structures}
\label{sec:03}

A generalized Riemannian metric $G$ is equivalent to a choice of a pseudo-Riemannian metric $g$ and a 2-form $F$ (a (1,1)-tensor $A$ satisfying \eqref{m1}), such that $G=g+F$.
A generalized metric connection, i.e. a linear connection preserving $G$, is a metric connection preserving the fundamental 2-form $F$. This supplies a number of examples (given below using weak metric structures defined by V.\,Rovenski and R.\,Wolak in \cite{rov-108}) with a (1,1)-tensor $A$ of constant rank.
 For $\nabla Q=0$ and the Nijenhuis (0,3)-tensor $N_Q(X,Y,Z):=g(N_Q(X,Y),Z)$, we~get
\begin{equation*}%\label{Qnuj2}
 N_Q(X,Y,Z){=} {-}T(QX,QY,Z){-}T(X,Y,Q^2Z){+}T(QX,Y,QZ){+}T(X,QY,QZ).
\end{equation*}
%}

\begin{definition}\rm
Let $(M, G = g + F)$ be a generalized Riemannian manifold.
A~linear connection \( \nabla \) is said to have {\it $Q$-torsion condition} if its torsion
%tensor
\( T \) satisfies
\begin{equation}\label{Eq-Q-T}
T(QX, Y) = T(X, QY)=QT(X,Y) \quad \text{for all } X, Y \in \mathfrak{X}_M,
\end{equation}
where \( Q: TM\to TM \) is an endomorphism that is self-adjoint with respect to the symmetric part $g$ of metric $G$, i.e., $g(QX, Y) = g(X, QY)$ for all $X,Y\in\mathfrak{X}_M$.
\end{definition}

Note that the $Q$-torsion condition is trivial when $Q=\lambda\,{\rm Id}$ for some real $\lambda>0$ (or in the classical case).
The $A$-torsion condition is not trivial when~$Q=\lambda\,{\rm Id}$.

{%\color{blue}
\begin{remark}\rm
If the torsion (0,3)-tensor is totally skew-symmetric, \eqref{st}, then the $Q$-torsion condition \eqref{Eq-Q-T} follows from the $A$-torsion condition \eqref{Eq-A-T}, but the converse does not hold in general.
\end{remark}
}

\subsection{Weak almost Hermitian structure} \label{sec:2.5.2}

Let us consider a weak almost Hermitian manifold $(M, g, A, Q)$, see \eqref{WAH}.

\begin{definition}\rm
%\textit{Weak nearly K\"ahler manifolds} are defined by
%a constraint only on the symmetric part of $A$: %$(\nabla^g_X A)X=0$ for $X\in\mathfrak{X}_M$.
%In other words,
A weak almost Hermitian manifold is said to be \textit{weak nearly K\"ahler} if the
covariant derivative of $A$ with respect to the Levi-Civita connection $\nabla^g$
%of the metric $g$
is skew-symmetric:
$$(\nabla^g_X A)X=0\quad \Longleftrightarrow \quad (\nabla^g_X F)(X,Y)=0.$$
If $\nabla^g A=0$, then such $(M, g, A, Q)$ is called a \textit{weak K\"ahler manifold}.
\end{definition}

At the same time, $(M, g, A, Q)$ admits a generalized Riemannian structure $G=g+F$.
Suppose that a linear connection \( \nabla \) on
%a generalized Riemannian manifold
$(M,G=g+F)$ preserves $G$, then it also preserves the metric \( g \), the form \(F\), and the tensor \(A\), see Theorem~\ref{main},
and we have \( \nabla g = \nabla F = \nabla A = 0 \). Consequently,
%the connection
\( \nabla \) also preserves the tensor \( Q \), i.e., \( \nabla Q = 0 \).
Thus, using the contorsion $K$, \eqref{skct}, and its relation with the torsion tensor,
$K(X,Y,Z):=g(K(X,Y),Z) = \frac{1}{2} T(X,Y,Z)$, we get
\begin{align}\label{E-nabla-g-A}
\notag
 (\nabla^g_X A)Y & = \nabla^g_X (AY) {-} A\nabla^g_X Y
 = A K(X,Y) {-} K(X,AY) \\
 & = \dfrac 12 \big[A T(X,Y) {-} T(X,AY)\big], \\
 \label{E-nabla-g-Q}
 (\nabla^g_X Q)Y & = Q K(X,Y) - K(X,QY) = \dfrac 12\big[Q T(X,Y) - T(X,QY) \big] .
\end{align}
For the case of totally skew-symmetric torsion (0,3)-tensor of $\nabla$, the $A$-torsion condition \eqref{Eq-A-T} gives
\begin{align}\label{QAQ}
% \begin{split}
 T(X,AY) & =T(AX,Y)=-AT(X,Y),
% \\  & T(QX,Y) =T(X,QY)=QT(X,Y),
%\end{split}
\\
%\end{equation}
%and
%\begin{equation}
\label{TQQ}
\notag
 T(X,QY) & =-T(AX,AY)=AT(X,AY) \\
 & =AT(AX,Y)=T(QX,Y)=QT(X,Y),
\end{align}
hence the $Q$-torsion condition \eqref{Eq-Q-T}.
%{%\color{red}\% MZ   have written a paragraph above about the contorsion tensor, so this might be redundant.\\
%Note that from the assumption $(\nabla_X\,g)(Y,Z)=0$ we get \eqref{st}:
%\begin{align*}
% K(X,Y,Z) + K(X,Z,Y)
 %= \frac{1}{2}\big[T(X,Y,Z)+T(X,Z,Y)\big] = 0.
%\end{align*}
%}
By (\ref{E-nabla-g-A}),
% we conclude that
%the torsion satisfying
the $A$-torsion condition implies
\begin{equation}\label{Fcond}
    (\nabla^g_{X}A)Y = -A T(X,Y).
\end{equation}
This leads us to the following theorem.

%Using (\ref{ATors}) and (\ref{WAH}), we obtain
%\begin{equation}\label{Eq-T-A}
%\begin{split}
%T(QX,AY,Z)&=T(AX,QY,Z)=T(QX,Y,AZ)=T(AX,Y,QZ)\\
%&=T(X,AY,QZ)=T(X,QY, AZ)=-T(AX,AY,AZ).
%\end{split}
%\end{equation}

\begin{theorem}\label{mainweak}
Let $(M, g, A, Q)$ be a weak almost Hermitian manifold with a fundamental 2-form \(F\), considered as a generalized Riemannian manifold \((M, G = g + F)\).
Suppose that \(\nabla\) is a generalized metric connection $($i.e., \(\nabla G = 0)\) with totally skew-symmetric torsion (0,3)-tensor $T$.
Then \(\nabla\) satisfies the $A$-torsion condition \eqref{Eq-A-T} if and only if \((M, g, A, Q)\) is a weak nearly K\"ahler manifold.
In~the case of \eqref{Eq-A-T}, the torsion
%tensor
of $\nabla$ is determined by the structure tensors as follows:
\begin{equation}\label{TTT}
\begin{split}
T(AX, Y, Z) &= -\dfrac 13\,dF(X, Y, Z),\\
T(QX, Y, Z) &= \dfrac 14\,N_A(X, Y, Z)
\ \big( =\dfrac 13\,dF(AX, Y, Z)\ \big).
\end{split}
\end{equation}
\end{theorem}

\begin{proof}
Assume that the connection \(\nabla\) with totally skew-symmetric torsion (0,3)-tensor satisfies the \(A\)-torsion condition \eqref{Eq-A-T}.
%Then, we have
%\begin{equation}\label{AAAQ}
% T(AX,Y,Z)=T(X,AY,Z)=T(X,Y,AZ).
%\end{equation}
%By applying  \eqref{conn2}, we obtain
%\begin{equation}\label{Feq}
% \begin{split}
% (\nabla^g_X F)(Y,Z) & + (\nabla^g_Y F)(X,Z) \\
%&= -\dfrac{1}{2} \big[ T(X,Y,AZ) + T(Z,X,AY) \big]
%   - \dfrac{1}{2} \big[ T(Y,X,AZ) + T(Z,Y,AX) \big] \\
%&= -\dfrac{1}{2} \big[ T(Z,X,AY) + T(Z,Y,AX) \big] \\
%&= -\dfrac{1}{2} \big[ T(X,AY,Z) - T(AX,Y,Z) \big].
%\end{split}
%\end{equation}
%From \eqref{AAAQ}, we obtain \((\nabla^g_X F)(Y,Z) + %(\nabla^g_Y F)(X,Z) = 0\),
From \eqref{Fcond} with $Y=X$ we conclude that
%and consequently,
\((M, g, A, Q)\) is a weak nearly K\"ahler manifold.
%{%\color{blue} \% VR: It seems that the "weak nearly %K\"ahler" property follows directly from \eqref{Fcond} %with $Y=X$, if so then we do not need \eqref{AAAQ} and %\eqref{Feq}.}
Conversely, assuming $(\nabla^g_XF)(X,Z)=0$ and taking into account (\ref{conn2}), we obtain
$T(AX,Y,Z)=T(X,AY,Z)$,
which implies that the linear connection $\nabla$ satisfies the $A$-torsion condition.
%\smallskip
Further, using \eqref{QAQ} in (\ref{conn1}), we obtain
\begin{equation*}%\label{TA}
T(AX,Y,Z)=-\dfrac 13\,dF(X,Y,Z).
\end{equation*}
Using (\ref{TQQ}) in (\ref{nuj2}) and (\ref{tor1}), we find that
$N_A$
%the Nijenhuis tensor of $A$
is totally skew-symmetric and
\begin{equation*}%\label{TQ}
T(QX,Y,Z)=\dfrac 14N_A(X, Y, Z)
\end{equation*}
is true,
that completes the proof of \eqref{TTT}.
\end{proof}

\begin{theorem}\label{Th-3.4}
Let the conditions of Theorem~\ref{mainweak} be satisfied.
%Let $(M, g, A, Q)$ be a weak almost Hermitian manifold %with a fundamental 2-form \(F\), considered as a %generalized Riemannian manifold \((M, G = g + F)\).
%Suppose that a linear connection \(\nabla\) with totally %skew-symmetric torsion tensor $T$ preserves the %generalized Riemannian metric \(G\), i.e., \(\nabla G = 0\).
Then the $Q$-torsion condition \eqref{Eq-Q-T}

\noindent\ \
(i) is equivalent to the following:
%the connection \(\nabla\) satisfies
%the $Q$-torsion condition \eqref{Eq-Q-T} is true if and %only if
the Levi-Civita connection preserves the tensor \(Q\), i.e., \(\nabla^g\,Q = 0\).

\noindent\ \
(ii) provides the following:
%the torsion tensor \(T\) satisfies the {\(Q\)-torsion %condition} if and only if
the Nijenhuis tensor \(N_Q\), see \eqref{nuj} with $P=Q$, vanishes, thus \(Q\) is integrable.
\end{theorem}

\begin{proof}
(i) Using~\eqref{E-nabla-g-Q}, we have
%\[
 $2\,(\nabla^g_X Q)(Y) = Q T(X,Y) - T(X, QY)$.
%\]
Therefore, \(\nabla^g Q = 0\) if and only if \(Q T(X,Y) = T(X, QY)\), which is exactly the $Q$-torsion condition in the case of totally skew-symmetric (0,3)-tensor $T$ and self-adjoint~$Q$:
\begin{align*}
 g(QT(X,Y),Z) & = g(T(X,Y), QZ) = T(X,Y,QZ) = T(Y,QZ, X) \\
 & = g(T(Y, QZ), X) = g(T(QY, Z), X) = T(QY, Z, X) \\
 & = T(X, QY, Z) = g( T(X,QY), Z).
\end{align*}
This completes the proof of (i).

(ii) Let the totally skew-symmetric torsion (0,3)-tensor satisfies the \(Q\)-torsion condition. So, we have
%\[
 $T(QX, Y, Z) = T(X, QY, Z) = T(X, Y, QZ)$.
%\]
Then
\begin{align*}
T(QX, QY, Z) & = T(X, Q^2 Y, Z) = T(X, Y, Q^2 Z), \\
T(QX, Y, QZ) & = T(X, QY, QZ) = T(X, Y, Q^2 Z).
\end{align*}
Substituting the last two equalities into \(N_Q\), we have
\begin{align*}
 N_Q(X, Y, Z) & = - T(X, Y, Q^2 Z) {-} T(X, Y, Q^2 Z) {+} T(X, Y, Q^2 Z) {+} T(X, Y, Q^2 Z) \\
 & = 0.
\end{align*}
%Conversely, if \(N_Q(X,Y,Z)=N_Q(X,Z,Y)=0\), then from %\eqref{Qnuj2} and (\ref{QAQ}), we obtain
%$T(Y, Q^2Z)=T(Z, Q^2Y)$ -- $Q^2$-torsion condition.
%Hence $$?$$
 This completes the proof of (ii).
\end{proof}

%What we can conclude from $N_Q=0$:
%\begin{itemize}
%  \item The operator \( Q \) is locally diagonalizable with real eigenvalues. \ \%VR: but not globally!
 % \item The eigen-distributions \( E_\lambda := \ker(Q - \lambda {\rm Id}) \) are smooth and involutive.
 % \item For each eigenvalue \( \lambda \), there exists a local integral submanifold \( \Sigma_\lambda \subset M \) such that \( T_p \Sigma_\lambda = E_\lambda(p) \).
 % \ \%VR: $\dim E_\lambda(p) = \mbox{multiplicity of } \lambda$ and $\Sigma_\lambda$ exists "globally".
 % \item In local coordinates adapted to these distributions, both the metric \( g \) and the tensor \( Q \) have diagonal form. \ \%VR: Not sure about coordinates. Is this the same as the 1-st item?
%\end{itemize}

%In other words, \( Q \) defines a local orthogonal decomposition of the manifold \( M \) into integral leaves of the eigendistributions, which are mutually orthogonal with respect to \( g \).

%\ \%VR: Fine. This is presented in Theorem 3.6!

\begin{example}\rm
The \emph{Chern connection} $\nabla^C$, see \cite{Chern}, referred to in some papers as the \emph{canonical} or \emph{Hermitian} connection is the unique linear connection on a Hermitian manifold $(M, g, A)$ preserving the generalized Riemannian metric $G = g+F$,
%The~fundamental form $F$ of the almost Hermitian %structure $(g,A)$ is defined by $F(X,Y)=g(AX,Y)$.
%We recall that the Chern connection $\nabla^C$ is the unique linear connection preserving the symmetric part
%of metric $g$ and the complex structure $A$,
so that the torsion tensor \(T\) of \(\nabla^C\) has the property \(T(X, AY) = T(AX, Y)\) for all $X,Y\in\mathfrak{X}_M$. This, implies
\[
T(AX, Y, Z) = T(X, AY, Z)=-T(X, Y,AZ) \quad \text{for all } X, Y, Z \in \mathfrak{X}_M .
\]
The connection $\nabla^C$ does not have totally skew-symmetric torsion (0,3)-tensor, and its torsion tensor satisfies different properties than the \(A\)-torsion condition (with an opposite sign in the term \(AT(X,Y)\) compared to \eqref{QAQ}).
This provides motivation for introducing the
$A$-torsion condition \eqref{Eq-A-T} and extending
%the Chern connection
$\nabla^C$ to the weak Hermitian structure $(g,A)$ on $M$.
Using $A^2=-Q$, we~get
\[
T(AX, AY, Z) = -T(QX, Y, Z)=-T(X, QY, Z)=-T(X, Y,QZ)
\]
for all $X, Y, Z \in \mathfrak{X}_M$.
The torsion tensor of $\nabla^C$ is determined by
\[
 T(QX,Y,Z) = -\dfrac 12\big[ dF(AX,Y,Z)+dF(X,AY,Z)\big]
\]
and the following is true:
%\[
$T(AX,Y,Z)=(\nabla^g_Z F)(X,Y)$.
%\]
\end{example}

\begin{example}\rm
The \emph{Bismut connection} $\nabla^B$, see \cite{Bismut}, also known as Strominger or Strominger-Bismut connection, is the unique linear connection
with totally skew-symmetric torsion (0,3)-tensor on a Hermitian ma\-nifold $(M, g, A)$ preserving the generalized Riemannian metric $G = g+F$.
%The fundamental form $F$ of the Hermitian structure %$(g,A)$ is defined by $F(X,Y)=g(AX,Y)$.
This provides motivation for extending
%the connection
$\nabla^B$ to the weak Hermitian structure $(g,A)$.
%on $M$.
Using $A^2=-Q$ and $N_A=0$, we have
%the following:
\[\begin{split}
 T(QX, Y, Z)&=T(X, QY, Z)=T(X, Y,QZ)\\
 &=T(AX,AY,Z)+T(AX,Y,AZ)+T(X,AY,AZ)
\end{split}
\]
for all $X, Y, Z \in \mathfrak{X}_M$.
From (\ref{tor1}), the torsion (0,3)-tensor of $\nabla^B$ is determined~by
\[
 T(AX,AY,Z) = dF(X,Y,AZ)
\]
and the following is true:
%\[
$dF(QX,Y,Z)=dF(X,QY,Z)=dF(X,Y,QZ)$.
%\]
%Note that,
In~general, on an almost weak Hermitian manifold, the connection $\nabla^B$ has no totally skew-symmetric torsion (0,3)-tensor.
\end{example}

The results obtained here provide a natural generalization of the previously established results for almost Hermitian (and almost para-Hermitian) manifolds \cite{FI,GKMW}. Our approach thus extends and unifies these results within more general ``weak" geometric structures.

The 6-dimensional unit sphere $\mathbb{S}^6$ in the set of
all purely imaginary Cayley numbers is a classical example of a non-K\"{a}hler nearly K\"{a}hler manifold. Other examples of 6-dimensional non-K\"{a}hler
nearly K\"{a}hler manifolds include $\mathbb{S}^3\times\mathbb{S}^3$,
the complex projective space $\mathbb{CP}^3$ and the flag manifold~$\mathbb{F}^3$.

The following example presents many
non-K\"{a}hler weak K\"{a}hler manifolds
as well as non-nearly K\"{a}hler weak nearly K\"{a}hler manifolds.

\begin{example}%\label{Ex-sHK}
\rm
L.\,P. Eisenhart \cite{E-1923} proved that if a Riemannian manifold $(M, g)$ admits a parallel symmetric 2-tensor other than the constant multiple of $g$, then it is reducible.
Therefore, a weak K\"{a}hler~manifold $(M, g, A, Q)$ with $Q\ne\lambda\,{\rm Id}$ for some $\lambda\in C^\infty(M)$ is reducible.
Take two or more almost Hermitian manifolds $(M_j,g_j, J_j)$, thus $J_j^2=-{\rm Id}_{\,j}$.
The~pro\-duct $\prod_{j}(M_j,g_j,\sqrt{\lambda_j}J_j)$, where $\lambda_j\ne1$ are
different positive constants, is a weak almost Hermitian manifold with $Q=\bigoplus_{j}\lambda_j{\rm Id}_{j}$.
We call $\prod_{\,j}(M_j,g_j,\sqrt{\lambda_j}J_j)$ a $(\lambda_1,\ldots,\lambda_k)$-\textit{weighed product of almost Hermitian manifolds};
in this case, $(M,g=\bigoplus_{j} g_j)$ decomposes into $k$ weak almost Hermitian~factors.
%%

%Moreover,
If $(M_j, g_j, A_j)$ are (nearly) K\"{a}hler manifolds, then their $(\lambda_1,\ldots,\lambda_k)$-weighed pro\-duct is a non-(nearly) K\"{a}hler weak (nearly) K\"{a}hler manifold.
By the above, any non-K\"{a}hler weak K\"{a}hler manifold is locally a weighed product of K\"{a}hler manifolds.
The classification of weak nearly K\"{a}hler manifolds in dimensions $\ge 4$ is an open problem.
%%%%%
Some 4-dimensional non-K\"{a}hler weak K\"{a}hler ma\-nifolds appear as $(\lambda_1,\lambda_2)$-weighed pro\-ducts of 2-dimensional K\"{a}hler manifolds.
Some 6-dimensional non-K\"{a}hler weak K\"{a}hler manifolds appear as $(\lambda_1,\lambda_2,\lambda_3)$-weighed products of 2-dimensional K\"{a}hler manifolds or $(\lambda_1,\lambda_2)$-weighed products of
2-dimensional and 4-dimensional K\"{a}hler manifolds, and similarly for even dimensions $>6$.
Some 8-dimensional non-nearly K\"{a}hler weak nearly
K\"{a}hler manifolds appear as $(\lambda_1,\lambda_2)$-weighed pro\-ducts of 2-dimensional K\"{a}hler manifolds and 6-dimensional nearly K\"{a}hler manifolds, and similarly for even dimensions~$>8$.

Such $(\lambda_1,\ldots,\lambda_k)$-weighed products of almost Hermitian manifolds can serve as new models of generalized Riemannian manifolds for theoretical physics.
\end{example}

\begin{theorem}\label{Th-001}
Let the conditions of Theorem~\ref{mainweak} be satisfied.

(i)~Suppose that $Q = \lambda\,{\rm Id}$ for some $\lambda\in C^\infty(M)$. Then $\lambda=const>0$,
the $Q$-torsion condition \eqref{Eq-Q-T} is true, and
$(\lambda^{-1/2}A,g)$ is an almost Hermitian structure; moreover, if the linear connection \(\nabla\) satisfies the $A$-torsion condition \eqref{Eq-A-T}, then
$(\lambda^{-1/2}A,g)$ is a nearly K\"{a}hler structure.

(ii)~Suppose that $Q\ne\lambda\,{\rm Id}$ for any $\lambda\in C^\infty(M)$. Then there exist $k>1$ mutually orthogonal even-dimensional distributions
${\mathcal D}_i\subset TM\ (1\le i\le k)$ such that $\bigoplus_{\,i}{\mathcal D}_i=TM$ and ${\mathcal D}_i$ are the eigen-distributions of $Q$
%(and $A$-invariant)
with different constant eigenvalues $\lambda_i:$  $0<\lambda_1<\ldots<\lambda_k$; moreover,
if the $Q$-torsion condition \eqref{Eq-Q-T}
$($or the $A$-torsion condition \eqref{Eq-A-T}$)$
is true, then the distributions ${\mathcal D}_i$ are involutive and define $\nabla^g$-totally geodesic foliations ${\mathcal F}_i$ and $M^{2m+s}(A,Q,g)$ is locally the $(\lambda_1,\ldots,\lambda_k)$-weighed product of almost Hermitian manifolds $($nearly K\"{a}hler manifolds, respectively$)$.
\end{theorem}

\begin{proof}
(i) Since $\nabla Q=0$, we get $\lambda=const>0$, hence
$(\lambda^{-1/2}A,g)$ is an almost Hermitian structure.
If the $A$-torsion condition
%\eqref{Eq-A-T}
is true,
then by Theorem~\ref{mainweak}, $(\lambda^{-1/2}A,g)$
is a nearly K\"{a}hler structure.

(ii) Suppose that $Q=-A^2$ is not conformal.
For a weak almost Hermitian structure $(g,A)$ on $M^{2m}$
we will build an $A$-\textit{basis} at a point $x\in M$
consisting of mutually orthogonal nonzero vectors of $T_xM$. Let $e_1\in T_xM$ be a unit eigenvector of the self-adjoint operator $Q>0$ with the minimal eigenvalue $\lambda_1\ne0$. Then, $Ae_1\in T_xM$ is orthogonal to $e_1$ and $Q(Ae_1) = A(Qe_1) = \lambda_1 Ae_1$.
Thus, the subspace of $T_xM$ orthogonal to the plane Span$\{e_1,Ae_1\}$ is $Q$-invariant (and $A$-invariant).
Note that $g(e_1, e_1)=g(Ae_1, Ae_1)
=g(Qe_1, e_1)=\lambda_1$.
Continuing in the same manner, we find a basis $\{e_1, Ae_1,\ldots, e_m, Ae_m\}$ of $T_x M$
consisting of mutually orthogonal vectors.
Hence, $Q$ has $k$ different eigenvalues $0<\lambda_1<\ldots<\lambda_k$ of even multiplicities $n_1,\ldots,n_k$,
and $\sum\nolimits_{\,i=1}^k n_i=2m$.
In this basis, $A$ and $Q$ have block-diagonal structures:
$Q = [\lambda_1{\rm Id}_{\,n_1},\ldots,\lambda_k{\rm Id}_{\,n_k}]$ and $A=[\sqrt{\lambda_1}J_{\,n_1},\ldots,\sqrt{\lambda_k}J_{\,n_k}]$, where $J_{\,n_i}$ is a complex structure on a $n_i$-dimensional subspace of~$T_xM$.
Since $Q$ is $\nabla$-parallel, we get the same structure at every point of $M$, i.e., $k$ and all $\lambda_i$ are constant on $M$, and
there exist mutually orthogonal $\nabla$-parallel (and $A$-invariant) eigen-distributions ${\mathcal D}_i$ of $Q$ with constant different eigenvalues $\lambda_i$.
Let the $Q$-torsion condition be true,
%Since $A$ is non-degenerate, from \eqref{QAQ}
%we get
%\begin{align}\label{Eq-T-A3}
% T(QX,Y)=T(X,QY)=Q\,T(X,Y).
%\end{align}
then using \eqref{Eq-Q-T} for any vector fields $X,Y\in{\mathcal D}_i$, i.e., $QX=\lambda_iX$ and $QY=\lambda_iY$, we have
%the following:
\begin{align*}
 Q[X,Y] & = Q\{\nabla_XY-\nabla_YX-T(X,Y)\}
 = \nabla_X(QY)-\nabla_Y(QX)-QT(X,Y) \\
 & =\lambda_i\{\nabla_XY-\nabla_YX\}-\lambda_i T(X,Y)
 =\lambda_i[X,Y].
\end{align*}
Hence each distribution ${\mathcal D}_i$ is involutive
and defines a foliation ${\mathcal F}_i$.
Similarly we can show that $Q(\nabla_XY)=\nabla_X(QY)=\lambda_i\nabla_XY$,
that is, each ${\mathcal F}_i$ is a $\nabla$-totally geodesic foliation.
%If the $Q$-torsion condition \eqref{Eq-Q-T} is true,
Then (by Theorem~\ref{Th-3.4}) $\nabla^g Q=0$,
hence
%${\mathcal D}_i$ is $T$-invariant; hence using the above
%and $\nabla_XY-%\nabla^g_XY=\frac12 T(X,Y)$, we conclude %that
each ${\mathcal F}_i$ is a $\nabla^g$-totally geodesic foliation and by de Rham Decomposition Theorem (see \cite{KN}), $(M,g)$ splits into the
$(\lambda_1,\ldots,\lambda_k)$-weighed product of
almost Hermitian manifolds.
If~the $A$-torsion condition \eqref{Eq-A-T} is true, then
by Theorem~\ref{mainweak}, $(M,g)$ splits into the $(\lambda_1,\ldots,\lambda_k)$-weighed product of nearly K\"{a}hler manifolds.\end{proof}

\subsection{Weak almost contact metric structure}

Contact Riemannian geomet\-ry is of growing interest due to its important role in both theoretical physics and pure mathematics. Weak almost contact structures, i.e., the complex structure on the contact distribution is approximated by a non-singular skew-symmetric tensor, allowed us to take a new look at the theory of contact manifolds and find new applications.
Let us consider a weak contact metric manifold $(M^{2m+1}, A, Q, \eta,\xi, g)$.
Denote by ${\mathcal D}=A(TM)$ a contact distribution on $M$.
Since $\nabla A=0$, the distribution
${\mathcal D}$ is $\nabla$-parallel; hence its unit normal, $\xi$, is also $\nabla$-parallel: $\nabla\xi=0$.
By duality,
%$\eta(\xi)=1$ and
$\eta=g(\xi,\cdot)$, we conclude that the 1-form $\eta$ is also $\nabla$-parallel: $\nabla\eta=0$.

On the other side, taking $\nabla$-derivative of $A\xi=0$ and using $\nabla A=0$, we obtain $A\nabla_X\,\xi=0$.
Hence $\nabla_X\,\xi\,||\,\xi$ for any $X$.
On the other hand, taking derivative of $g(\xi,\xi)=1$ and using $\nabla g=0$, we obtain $g(\nabla_X\xi,\xi)=0$. Therefore, $\nabla_X\xi=0$ for any $X$, i.e., $\nabla\xi=0$.
Then \eqref{acon} (see also Theorem~\ref{main}) gives $\nabla Q=0$:
\begin{align*}%\label{acon1}
%\notag
0 &=(\nabla_XA)AY + A(\nabla_XA)Y = (\nabla_X A^2)Y \\
&=-(\nabla_X Q)Y +(\nabla_X\,\eta)(Y)\,\xi
-\eta(Y) \nabla_X\,\xi =- (\nabla_X Q)Y.
\end{align*}
% conclude that Q is also parallel
Since $T$ is totally skew-symmetric, the equation (\ref{skct}) yields
$0 = g(\nabla^g_X\,\xi, Y) + \dfrac12 T(X, \xi, Y),$ which shows that $\xi$ is a geodesic and Killing vector field:
$\nabla^g_\xi\,\xi=0$ and
$g(\nabla^g_X\,\xi, Y)+g(\nabla^g_Y\,\xi, X)=0$.
Further, we have
\begin{align}\label{wdeta}
\notag
 d\eta(X,Y) & = (\nabla^g_X\,\eta)(Y) -(\nabla^g_Y\,\eta)X) \\
 & =g(\nabla^g_X\,\xi, Y)-g(\nabla^g_Y\,\xi, X)=T(X,Y,\xi)
\end{align}
and
\[
 d\eta(X,\xi)=-d\eta(\xi,X)=0.
\]
Now, using (\ref{acon}) in equation (\ref{nuj2}), the Nijenhuis tensor becomes:
\begin{align}\label{wNuj}
\notag
 N_A(X, Y, Z) + \eta(Z)d\eta(X, Y ) & = T(X, Y, QZ) -T(AX, AY, Z) \\
 & -T(AX, Y, AZ) -T(X, AY, AZ).
\end{align}
The right-hand side of \eqref{wNuj} is totally skew-symmetric if and only if the torsion tensor satisfies the $Q$-torsion condition \eqref{Eq-Q-T}:
%\[
$T(X, Y, QZ) = T(X, QY, Z)$.
%\]
Let us introduce the (0,3)-tensor
\begin{equation}\label{nwac}
N_A^{{\rm wac}} = N_A + d\eta \otimes \eta,
\end{equation}
which we will call the \emph{Nijenhuis tensor in the weak almost contact geometry}, similarly to the definition given by Blair (see, e.g., \cite{Blair}, where this formula is with the coefficient 2) in the case of almost contact structures.
%{%\color{blue} \% VR: In \cite{Blair}, p.~81, a similar formula has the factor "2":
%$N_A^{{\rm wac}}(X,Y) = N_A(X,Y) + 2\,d\eta(AX,Y)\xi$.
%%Another question: $\otimes$ or $\wedge$ (wedge product of forms)?
%%Perhaps when we move to (0,3)-tensors, the factor "2"  disappears? We can explain this in the paper.
%\\
%\% In our (3.18), Blair used 1/2 in $
%d\eta(X,Y) = 1/2((\nabla^g_X \eta)(Y) - (\nabla^g_Y \eta)(X))
%$ on p.~62 of the second edition of Blair's book. I %believe it is not necessary to elaborate further, as %this is purely a matter of notation.
%It should be $\otimes$. The symbol $\wedge$ represents %the antisymmetric wedge product of differential forms.
%}
From (\ref{wNuj}) and (\ref{wdeta}) and taking into account (\ref{acon}), we obtain
\begin{equation*}
 %\begin{split}
   N_A(X,Y,\xi)=-d\eta(AX,AY), \quad
   N_A(\xi, Y, Z) = d\eta(Y,QZ)-d\eta(AY, AZ).
 %\end{split}
\end{equation*}
So, the last equations give
\begin{equation*}
    \begin{split}
        d\eta(Y,QZ)=N_A(\xi, Y, Z) -N_A(Y,Z,\xi).
    \end{split}
\end{equation*}

The above considerations yield the following theorem.

\begin{theorem}\label{T-2.13}
Let $(M^{2m+1}, A, Q, \eta,\xi, g)$ be a weak almost contact metric manifold considered as a generalized Riemannian manifold \((M, G = g + F)\). Suppose that \(\nabla\) is a generalized metric connection $($i.e., \(\nabla G = 0)\) with totally skew-symmetric torsion $(0,3)$-tensor $T$. Then the following properties hold:
\begin{enumerate}
 \item[i)] The Reeb vector field $\xi$ is a geodesic and Killing vector field.
 \item[ii)] The connection $\nabla$ satisfies the $Q$-torsion condition~\eqref{Eq-Q-T} if and only if the Nijenhuis tensor in weak almost contact geometry $N^{{\rm wac}}$ given by ~\eqref{nwac} is totally skew-symmetric.
\end{enumerate}
\end{theorem}

%{%\color{red}\% VR: Good! A weak contact metric manifold (i.e., $d\eta=F$ - can we show this?) with a Killing $\xi$ is called {\em weak K-contact manifold} - an important class!

%\% MZ $d\eta=F$ is not valid in general for our case.}

The following theorem is similar to Theorem~\ref{Th-001}.
%is a particular case, $s=1$, of Theorem~\ref{Th-002}.

\begin{theorem}\label{Th-3.13}
Let the conditions of Theorem~\ref{T-2.13} be satisfied.
Then the following properties hold.

\noindent\ \
(i)~If $Q|_{\mathcal D}=\lambda\,{\rm Id}_{\mathcal D}$ for some $\lambda\in C^\infty(M)$, then $\lambda=const$ and $(\lambda^{-1/2}A,\eta,\xi,g)$ is an almost contact metric structure on $M^{2m+1}$;
%and $M^{2m+1}(\lambda^{-1/2}A,\eta,\xi,g)$ is locally the product of $\mathbb{R}$ and an almost Hermitian manifold;
moreover, if the connection \(\nabla\) satisfies
the $Q$-torsion condition \eqref{Eq-Q-T} $($or the $A$-torsion condition \eqref{Eq-A-T}$)$,
%the $A$-torsion condition \eqref{Eq-A-T},
then $(M^{2m+1}, \lambda^{-1/2}A,\eta,\xi,g)$ is locally the product of $\mathbb{R}$ and
an almost Hermitian manifold $($a nearly K\"{a}hler manifold, respectively$)$.
%a nearly K\"{a}hler manifold.

\noindent\ \
(ii)~If~$Q|_{\mathcal D}\ne\lambda\,{\rm Id}_{\mathcal D}$ for any $\lambda\in C^\infty(M)$, then there exist $k>1$ mutually orthogonal even-dimensional distributions
${\mathcal D}_i\subset{\mathcal D}$ such that $\bigoplus_{\,i=1}^k{\mathcal D}_i={\mathcal D}$ and ${\mathcal D}_i$ are eigen-distributions of $Q$ with constant eigenvalues $\lambda_i:$  $0<\lambda_1<\ldots<\lambda_k$;
moreover, if the $Q$-torsion condition \eqref{Eq-Q-T}
$($or the $A$-torsion condition \eqref{Eq-A-T}$)$
is true, then the distributions ${\mathcal D}_i$ are involutive and define $\nabla^g$-totally geodesic foliations ${\mathcal F}_i$
%and each distribution ${\mathcal D}_i$ is involutive and defines a $\nabla$-totally geodesic foliation ${\mathcal F}_i$.
%%
%Moreover, if the connection \(\nabla\) satisfies the $Q$-torsion condition \eqref{Eq-Q-T} $($or the $A$-torsion condition \eqref{Eq-A-T}$)$, then foliations ${\mathcal F}_i$ are $\nabla^g$-totally geodesic
and $(M^{2m+1}, A,Q,\eta,\xi,g)$ is locally a $(1,\lambda_1,\ldots,\lambda_k)$-weighed product of $\mathbb{R}$ and $k$ almost Hermitian manifolds $($nearly K\"{a}hler manifolds, respectively$)$.
\end{theorem}

%\begin{proof}
%This is similar to the proof of Theorem~\ref{Th-001}.
%\end{proof}

\begin{proof}
(i) From $\nabla Q=0$ we get $\lambda=const>0$ and $\lambda\ne1$, hence
$(\lambda^{-1/2}A, \eta,\xi, g)$ is an almost contact metric structure.
If the linear connection $\nabla$ satisfies the $Q$-torsion condition, then using \eqref{Eq-Q-T} for any vector fields $X,Y\in{\mathcal D}$, i.e., $QX=\lambda X$ and $QY=\lambda Y$ with $\lambda\ne1$, we have
%the following:
\begin{align*}
 Q[X,Y] & = Q\{\nabla_XY-\nabla_YX-T(X,Y)\}
 = \nabla_X(QY)-\nabla_Y(QX)-QT(X,Y) \\
 & =\lambda\{\nabla_XY-\nabla_YX\}-\lambda T(X,Y)
 =\lambda[X,Y].
\end{align*}
Hence the contact distribution ${\mathcal D}$ is involutive and defines a foliation ${\mathcal F}$.
Similarly we can show that $Q(\nabla_XY)=\nabla_X(QY)=\lambda\nabla_XY$,
that is, ${\mathcal F}$ is a $\nabla$-totally geodesic foliation.
Then (by Theorem~\ref{Th-3.4}) $\nabla^g Q=0$,
hence ${\mathcal F}$ is $\nabla^g$-totally geodesic %foliation
and by de Rham Decomposition Theorem, see \cite{KN},
$(M, \lambda^{-1/2}A, \eta,\xi, g)$ is locally the product of $\mathbb{R}$ and an almost Hermitian manifold.
If the linear con\-nection
$\nabla$ satisfies the $A$-torsion condition \eqref{Eq-A-T}, then
%by Theorem~\ref{Th-3.4},
%$M^{2m+1}(\lambda^{-1/2}A, \eta,\xi, g)$ is locally the %product of $\mathbb{R}$ and
the second factor is
a nearly K\"{a}hler~manifold.

(ii) Suppose that $Q|_{\mathcal D}$ is not conformal.
For a weak almost contact metric manifold, as well as for a weak almost Hermitian mani\-fold, see the proof of Theorem~\ref{Th-001},
there exists an $A$-\textit{basis} $\{e_1, Ae_1,\ldots, e_m, Ae_m; \xi\}$ of $T_x M$ consisting of mutually orthogonal nonzero vectors.
Hence, $Q$ restricted on ${\mathcal D}_x$ has $k>1$ different eigenvalues $0<\lambda_1<\ldots<\lambda_k$ of even multiplicities $n_1,\ldots,n_k$,
and $\sum\nolimits_{\,i=1}^k n_i= 2m$.
In this basis, $A$ and $Q$ restricted on ${\mathcal D}_x$ have block-diagonal structures:
$Q|_{\mathcal D} = [\lambda_1{\rm Id}_{n_1},\ldots,\lambda_k{\rm Id}_{n_k}]$ and $A|_{\mathcal D}= [\sqrt{\lambda_1}J_{n_1},\ldots,\sqrt{\lambda_k}J_{n_k}]$,
where $J_{n_i}$ is a complex structure on a $n_i$-dimensional subspace of~${\mathcal D}_x$.
Since $Q$ is $\nabla$-parallel, we get the same structure at each point of $M$, i.e., $k$ and all $\lambda_i$ are constant on $M$, and
there exist mutually orthogonal $\nabla$-parallel
(and $A$-invariant)
eigen-distributions ${\mathcal D}_i\subset{\mathcal D}$ of $Q$ with constant different eigenvalues $\lambda_i$.
The rest of the proof is similar to the proof of Theorem~\ref{Th-001}.
%\% CONT.
\end{proof}

\begin{remark}\rm
Condition $\nabla G=0$ (hence $\nabla g=\nabla A=0$) in Theorem~\ref{Th-3.13} is sufficiently strong.
When the contact distribution ${\mathcal D}$ has constant dimension $<\dim M$, it makes sense to examine the weaker condition: $\nabla G=0$ along~${\mathcal D}$.
\end{remark}

\begin{remark}[Weak $f$-structure]\rm
%\subsection{Weak \texorpdfstring{$f$}{f}-structure, %\texorpdfstring{$f^3+fQ=0$}{f extasciicircum 3+fQ=0}}

An~$f$-structure by K.\,Yano \cite{yan}, where $f^3+f=0$,
on a smooth manifold $M^{2m+s}$ serves as a higher-dimensional analog of {almost complex structures} ($s=0$) and {almost contact structures} ($s=1$).
The following generalizes the $f$-structure.
Let us consider a \textit{weak $f$-manifold} $(M, g, A, Q)$, i.e., a
%(pseudo-)
Riemannian manifold $(M, g)$ of dimension $n=2m+s$ endowed with a (1,1)-tensor $f=A$ of rank $2m$, and a non-singular (1,1)-tensor $Q$, satisfying the following compatibility conditions, see~\cite{rov-108,rst-139}:
\begin{equation}\label{Eq-f}
 A^3+AQ=0,\quad g(AX,A^2Y)=g(QX,AY),\quad F(X,Y)=g(AX,Y).
\end{equation}
Note that if $A$ is non-degenerate ($s=0$), then
$A^3+AQ=0$ reduces to $A^2=-Q$.
%, see Section~\ref{sec:2.5.2}.

By the first equality in \eqref{Eq-f}, $Q$ commutes with $A$: $[A,Q]=0$.
Due to \eqref{Eq-f}, we have complementary orthogonal distributions
${\mathcal D}=A(TM)$ and ${\mathcal D}_0=\ker A$,
where $\dim{\mathcal D}=2m$ and $\dim{\mathcal D}_0=s$.
We omit here our result with $s>1$,
for a weak $f$-manifold with a fundamental 2-form \(F\), considered as a generalized Riemannian manifold \((M^{2m+s}, G = g + F)\); it is similar to Theorem~\ref{Th-3.13}.
\end{remark}

%\% Sure. Very good consideration.
%\% VR: Obtaining some solutions to the above problems, we could finish this paper.
%\% Okay, our next topic is Einstein metricity condition (NGT)
%\% If you are interested? It would be good to continue with the idea.
%I agree. First, we will finish this paper (concentrate on basic questions), then think on the theme in general and grant proposal (with several topics on this subject).
%So, let us continue!
%Yes, let us finish the paper first. I will add some sentence in the introduction.
%Let's think about where to submit the paper. (Filomat, MDPI?) Arxiv is required.
%\% VR: MDPI is possible.

\subsection{Weak almost para-Hermitian structure and weak para-$f$-structure}
\label{sec:3.3}

A \textit{weak almost para-Hermitian manifold}, $(M,g,A,Q)$, is a pseudo-Riemanni\-an manifold $(M,g)$ of dimension $2m\ge4$ and signature $(m,m)$ endowed with non-singular endomorphisms: $A$ (skew-symmetric) and $Q\ne{\rm Id}$ (self-adjoint) and the fundamental 2-form $F$, such that the following conditions hold, see \cite{rov-137}:
\begin{equation*}%\label{APH}
A^2=Q, \quad g(AX,AY)=-g(QX,Y),\quad F(X,Y)=g(AX,Y) .
\end{equation*}
This structure generalizes the almost para-Hermitian structure, see \cite{CFG-1996}.

The results analogous to Theorem~\ref{mainweak}
%{\color{blue} [and \ref{Th-001} ?]}
from Section \ref{sec:2.5.2} also hold for the class of weak almost para-Hermitian manifolds, and their statements and proofs can be applied here accordingly.

The difference in signs appears in formula (\ref{TTT}), where, under the conditions of Theorem~\ref{mainweak}, the Nijenhuis tensor is totally skew-symmetric and satisfies
\begin{equation*}%\label{para-TTT}
T(QX, Y, Z) = -\dfrac 14 N_A(X, Y, Z)
\ \big( =-\dfrac 13 dF(AX, Y, Z)\ \big).
\end{equation*}
%%%%%%%%%%%%%%%%%%%%%%%%%%%
%\subsection{Weak para-$f$-structure, $f^3-fQ=0$}
The following generalizes the classical para-$f$-structure, where $f^3-f=0$.
A~weak para-$f$-manifold $(M, g, A, Q)$ is a pseudo-Riemannian manifold $(M, g)$ of dimension $2m+s$ endowed with a (1,1)-tensor $f=A$ of rank $2m$, and a non-singular (1,1)-tensor $Q\ne{\rm Id}$, satisfying the compatibility conditions, see~\cite{rov-121}:
\begin{equation}\label{Eq-f-para}
 A^3-AQ=0,\quad g(AX,A^2Y)=-g(QX,AY),\quad F(X,Y)=g(AX,Y).
\end{equation}
Note that if $A$ is non-degenerate ($s=0$), then
$A^3-AQ=0$ reduces to $A^2=Q$.
%, see Section~\ref{sec:2.5.3}.
Then  ${\mathcal D}=A(TM)$ and ${\mathcal D}_0=\ker A$ are two complementary orthogonal distributions on $M$.
We~have ${\mathcal D}={\mathcal D}^+\oplus{\mathcal D}^-$,
where ${\mathcal D}^+$ and ${\mathcal D}^-$ are positive and negative eigen-distributions of $Q$ restricted to ${\mathcal D}$,
with $\dim{\mathcal D}^+=\dim{\mathcal D}^-=m$ and $\dim{\mathcal D}_0=s$.
By the first equality in \eqref{Eq-f-para}, $Q$ commutes with $A$: $[A,Q]=0$.
%We omit the result similar to Theorem~\ref{Th-3.13}.
%\subsubsection
{Weak metric para-$f$-manifolds} are introduced and can be examined similarly to weak metric $f$-manifolds.

\smallskip

%\subsubsection
%\textbf{Weak almost para-contact metric structure}

We pay special attention to weak almost para-contact metric manifolds $(M^{2m+1}, A, Q, \eta,\xi, g)$,  i.e. a $(2m+1)$-dimensional pseudo-Riemannian manifold of signature $(m+1,m)$ equipped with a (1,1)-tensor $A$ of rank $2m$,
a 1-form $\eta$, a vector field $\xi$ dual to $\eta$ with respect to the metric $g$,\ $\eta(\xi)=1,\ \eta(X)=g(X,\xi)$, and a non-singular (1,1)-tensor $Q$, satisfying the following
%compatibility
conditions:
\begin{equation}\label{apcon}
A^2=Q-\eta\otimes\xi, \ \
g(AX,AY)=-g(QX,Y)+\eta(X)\,\eta(Y), \ \ A\xi=0, \ \ Q\,\xi=\xi.
\end{equation}
{%\color{blue}
For $Q={\rm Id}$, this defines almost para-contact metric manifolds, for example, \cite{Calv,IZ,Zam}.
}
In the case \eqref{apcon} the skew-symmetric part $F$ of $G=g+F$ has~rank equal to $2m$, and the equality $F(\xi,X)=0$ holds.
%\smallskip
In a similar way, for the totally skew-symmetric torsion (0,3)-tensor $T$ of a metric connection $\nabla$ on $M$, we conclude that $\xi$ is both a geodesic and a Killing vector field and $d\eta(X,Y)=T(X,Y,\xi)$ holds.
By applying (\ref{apcon}) in \eqref{nuj2}, the Nijenhuis tensor of~$A$ becomes:
\begin{align}\label{wpNuj}
\notag
 N_A(X, Y, Z)-\eta(Z)d\eta(X, Y ) & = -T(X, Y, QZ) -T(AX, AY, Z) \\
 & -T(AX, Y, AZ) -T(X, AY, AZ).
\end{align}
The right-hand side of \eqref{wpNuj} is totally skew-symmetric if and only if the torsion tensor of $\nabla$ satisfies the $Q$-torsion condition \eqref{Eq-Q-T}.
In the same manner, we introduce the (0,3)-tensor
\begin{equation}\label{nwapc}
N_A^{{\rm wapc}} = N_A -d\eta \otimes \eta,
\end{equation}
which we will call the \emph{Nijenhuis tensor in weak almost para-contact geometry}.

From (\ref{wpNuj}) and $d\eta(X,Y)=T(X,Y,\xi)$, taking into account (\ref{apcon}), we obtain
\begin{equation*}
 %\begin{split}
 N_A(X,Y,\xi)=-d\eta(AX,AY), \quad
 N_A(\xi, Y, Z) = d\eta(Y,QZ)-d\eta(AY, AZ).
% \end{split}
\end{equation*}

%Thus, the following theorems (similar to Theorems \ref{T-2.13} and \ref{Th-3.13}) hold.
Thus, the following theorem (similar to Theorem \ref{T-2.13}) is true.

\begin{theorem}%\label{T-para-2}
Let $(M^{2m+1}, A, Q, \eta,\xi, g)$ be a weak almost para-contact metric manifold considered as a generalized Riemannian manifold \((M, G = g + F)\).
Suppose that there exists a generalized metric connection \(\nabla\),
%preserving the generalized Riemannian metric \(G\),
i.e., \(\nabla G = 0\), with a totally skew-symmetric torsion $(0,3)$-tensor. Then the following~hold:
\begin{enumerate}
 \item[i)] the Reeb vector field $\xi$ is a geodesic and Killing vector field.
 \item[ii)] the connection $\nabla$ satisfies the $Q$-torsion condition~\eqref{Eq-Q-T} if and only if the Nijenhuis tensor in weak almost para-contact geometry $N^{{\rm wapc}}$ given by ~\eqref{nwapc} is totally skew-symmetric.
\end{enumerate}
\end{theorem}

{%\color{blue}
\begin{proof}
 i) Taking the covariant derivative of $A\xi = 0$ and using $\nabla A = 0$, we obtain $A \nabla_X \xi = 0.$
Hence, $\nabla_X \xi \parallel \xi$  for any $X$.
On the other hand, differentiating $g(\xi, \xi) = 1$ and using $\nabla g = 0$, we obtain
%\[
 $g(\nabla_X \xi, \xi) = 0$.
%\]
Hence,
 $\nabla_X \xi = 0$ for $X\in\mathfrak{X}_M$, i.e., $\nabla \xi = 0$.
Putting $X = \xi$ in the previous equality yields
%\[
$\nabla_\xi\,\xi = 0$,
%\]
and thus the Reeb vector field $\xi$ is a geodesic and Killing vector field.

ii) Assuming that $N_A^{\mathrm{wapc}}$ is totally skew-symmetric, we get
\[
T(X,Y,QZ) = T(X,QY,Z).
\]
Since $T$ is skew-symmetric in $X$ and $Y$, it satisfies the $Q$-equation \eqref{Eq-Q-T}.
Conversely, if the $Q$-equation \eqref{Eq-Q-T} holds, then
\[
\begin{aligned}
N_A^{\mathrm{wapc}}(X,Y,Z) &{=} {-}T(X,Y,QZ) {-} T(AX,AY,Z) {-} T(AX,Y,AZ) {-} T(X,AY,AZ) \\
&= \!T(X,Z,QY) {+} T(AX,Z,AY) {+} T(AX,AZ,Y) {+} T(X,AZ,AY) \\
&= -N_A^{\mathrm{wapc}}(X,Z,Y),
\end{aligned}
\]
which completes the proof.
\end{proof}
}

\begin{remark}\rm
Results similar to Theorems \ref{Th-001} and \ref{Th-3.13} cannot hold for a pseudo-Riemannian metric $g$, not a Riemannian one, since the spectral theorem is false unless $g$ is positive definite.
For any
%$n$-dimensional
pseudo-Riemannian manifold, results based on the diagonalizability of a self-adjoint operator $Q$ may be incorrect due to the possible existence of light-like eigenvectors.
In particular, the spectral theorem for self-adjoint operators is generally false unless the operator has no light-like eigenvectors, e.g. B.~O'Neill \cite[pp.~260-262]{ONeil}.
%and also K.~Rajaratnam~\cite{Raj} for a complete classification in Lorentzian signature.
%%
In a pseudo-Euclidean vector space $(V, g)$ of dimension $\ge 3$, a self-adjoint operator $Q: V \to V$ is diagonalizable with respect to a $g$-orthonormal basis only if
\(
g(QX, X) \ne 0
\)
for all light-like vectors $X \in V$, see \cite{Bognar}.
Thus, in the pseudo-Riemannian case, principal directions and principal curvatures are generally not defined without additional restrictions (e.g., spacelike or timelike hypersurfaces, where light-like vectors are excluded).
\end{remark}

%\begin{remark}\rm
%Over a pseudo-Riemannian manifold, the decomposition may require complexification due to the presence of lightlike eigenvectors; the distributions ${\mathcal D}_i$ correspond to generalized eigen-distributions of $Q$.
%Let $(M,g)$ be a pseudo-Riemannian vector space and $Q$ a self-adjoint endomorphism.
%Consider the operator $Q_\mathbb{C}$ associated with $Q$ (and defined on the complexification of $TM$).
%The~characteristic polynomial of $Q_\mathbb{C}$ splits over $\mathbb{C}$, and we obtain the decomposition
%\[
%$TM_\mathbb{C} = \bigoplus\nolimits_{\,\lambda, \Im(\lambda)\ge 0} (TM_\mathbb{C})_\lambda$,
%\]
%where $(TM_\mathbb{C})_\lambda$ denotes the generalized eigenspace corresponding to $\lambda \in \mathbb{C}$.
%Since the direct sum is orthogonal with respect to $g$, the generalized eigenspaces $TM_\lambda$ are nondegenerate,
%and each distribution $TM_\lambda$ can be used to define invariant subbundles of the tangent bundle for pseudo-Riemannian manifolds, possibly after %complexification.
%\end{remark}

%Analysis of theorem (similar to Theorem \ref{Th-3.13}) is rather complicated and we omit it here.

\section{Conclusion}

This paper presents new applications of weak contact metric structures to generalized Riemannian manifolds $(M,G=g+F)$ equipped with a metric connection of totally skew-symmetric torsion tensor.
We~described the behavior of the Nijenhuis tensor for a skew-symmetric non-singular endomorphism~$A$ and a symmetric non-singular endomorphism $Q$ with respect to symmetric part of metric in specific cases.
The vanishing of the Nijenhuis tensor $N_Q$, ensured by the $Q$-torsion condition, indicates that $Q$ defines a Nijenhuis operator. This naturally leads to the consideration of Nijenhuis manifolds, e.g. \cite{Matv}.
We~have shown that weighted products of almost Hermitian %and almost contact
manifolds can serve as new models of generalized Riemannian manifolds that can be useful for theoretical physics.
Following the idea, in the next paper we will consider weak structures on the generalized Riemannian manifolds satisfying the Einstein metricity condition, e.g.~\cite{IZ1}. %In addition, we will explore the case when $A$ is singular.
It would also be interesting
to explore the pseudo-Riemannian case and
to characterize generalized Riemannian manifolds equipped with weak metric structures of various distinguished classes, e.g. weak (almost) ${\mathcal K}$-, weak (almost) ${\mathcal C}$-,
weak (almost) ${\mathcal S}$-, and weak $f$-${\mathcal K}$-contact structures, see \cite{rst-139},
in terms of the totally skew-symmetric torsion tensor of a metric connection on $(M,G=g+F)$.

\bigskip

\noindent {\bf Acknowledgments}.
The authors express their sincere gratitude to the anonymous referees for a number of insightful comments, in particular, for pointing out the limitations of extending results similar to Theorems \ref{Th-001} and \ref{Th-3.13} to the pseudo-Riemannian case. These comments significantly improved both the accuracy and clarity of the manuscript.

M.Z. was partially supported by the Ministry of Education, Science and Technological Development of the Republic of Serbia,
project
No. 451-03-137/2025-03/200124 and by the Bulgarian Ministry of Education and Science, Scientific Programme ``Enhancing the Research Capacity in Mathematical Sciences (PIKOM)", No. DO1-67/05.05.2022.

\end{document}